\documentclass[12pt]{amsart}

\usepackage{amsmath}
\usepackage{amssymb} 
\newcommand{\Z}{{\mathbb Z}} 
\newcommand{\Q}{{\mathbb Q}} 

\newtheorem{thm}{Theorem}
\newtheorem*{ED} {Theorem (Edixhoven)}
\newtheorem*{JE} {Theorem (Jenkins)}
\newtheorem{prop}{Proposition}

\theoremstyle{remark}

\newcommand{\itop}[2]{\genfrac {}{}{0pt}{3}{#1}{#2} }

\begin{document}  

\author{P. Guerzhoy$^*$} \thanks{$^*$Partially supported by NSF grant DMS-0501225}
\address{Department of Mathematics, Temple University,
1805 N. Broad St., Philadelphia, PA 19122}
\email{pasha@math.temple.edu}
\title{Some congruences for traces of singular moduli}

\begin{abstract}
We address a question posed by Ono \cite[Problem 7.30]{Onobook}, prove a general
result for powers of an arbitrary prime, and provide an explanation for 
the appearance of higher congruence moduli  for certain small primes. 
One of our results overlaps but does not coincide with a recent
result of Jenkins \cite{Jenkins}. This result essentially coincides
with a recent result of Edixhoven \cite{Edixhoven},
and we hope that the comparison of the methods, which are entirely different, may
reveal a connection between the 
$p$-adic geometry and the arithmetic of half-integral weight Hecke operators.
\end{abstract}
\maketitle

\vspace{3mm}
\begin{center}{\sc 1. Introduction and discussion of the results} \end{center}
\vspace{3mm}

Throughout the paper 
$D$ and $d$ denote positive and non-negative integers, respectively,
which satisfy the congruences
$$
D \equiv 0,1 \bmod 4 \hspace{5mm} \textup{and} \hspace{5mm} d \equiv 0,3 \bmod 4.
$$
We denote by $\chi_d= \left( \frac{-d}{\cdot} \right)$
and $\chi_D = \left( \frac{D}{\cdot} \right)$ the  quadratic 
Dirichlet characters associated with 
the imaginary and real quadratic fields $\Q(\sqrt{-d})$ and $\Q(\sqrt{D})$ correspondingly.
The character $\chi_d$ (resp. $\chi_D$) is primitive if $-d$ (resp. $D$) is a fundamental
discriminant.

Zagier considered in \cite{Zagier} the nearly holomorphic modular forms 
$$
f_d=q^{-d} + \sum_{D>0}A(D,d)q^D
$$
of weight $1/2$ and
$$
g_D=q^{-D} + \sum_{d \geq 0}B(D,d)q^d.
$$
of weight $3/2$. (Here and in the following $q=\exp(2 \pi i \tau)$ with $\Im (\tau) >0$.)
The explicit recursive construction of $f_d$ and $g_D$, 
provided by Zagier, guarantees that $A(D,d),B(D,d) \in \Z$.

For an integer $m \geq 1$
the Hecke operators $T(m)$  act on these forms and conserve the
integrality of the Fourier coefficients. Following \cite{Zagier} we denote by
$A_m(D,d)$ and $B_m(D,d)$ the coefficient of $q^D$ in $f_d\vert_{\frac{1}{2}}T(m)$ 
and the coefficient of $q^D$ in $g_D\vert_{\frac{3}{2}}T(m)$, respectively.
Zagier proved that
\begin{equation} \label{AB}
A_m(D,d) = - B_m(D,d),
\end{equation}
and this common value divided by $\sqrt{D}$ is the  (twisted if $D>1$) trace of a certain modular
function. This interpretation in terms of traces of singular moduli is a primary source 
of motivation for the investigation of these numbers. 

Ahlgren and Ono  studied the arithmetic of traces of singular moduli in \cite{AO}, and, in particular,
proved the congruences
\begin{equation} \label{AO}
A_m(1,p^2d) \equiv 0 \bmod p
\end{equation}
if the prime $p$ splits in $\Q(\sqrt{-d})$ and $p \nmid m$. In this connection 
Ono posed a question \cite[Problem 7.30]{Onobook} whether there are natural generalizations of
(\ref{AO}) modulo arbitrary powers of $p$. Numerical evidence indicates that if
$\chi_d(p)=\chi_D(p)$, then 
\begin{equation} \label{target}
A_m(D,p^{2n}d) \equiv 0 \bmod p^n
\end{equation}
with the maximum congruence moduli which exceeds $p^n$ for $p \leq 11$.

The question splits into two parts:
to find similar congruences which hold for powers of an arbitrary prime and to find
series of stronger congruences for special primes. 

We firstly comment on the former part of the question. Recently Edixhoven \cite{Edixhoven} 
used the interpretation of the numbers
$A_m(1,p^{2n}d)$ as traces of singular moduli and the local moduli theory of 
ordinary elliptic curves  in positive characterstic and obtained the following 
result.

\begin{ED}
If $D=1$ and $\chi_d(p)=1$, then \textup{(\ref{target})} holds for any $m \geq 1$.
\end{ED}

Recently Jenkins \cite{Jenkins} presented an elementary argument based on the identity
(\ref{AB}) and standard formulas for the action of half-integral weight Hecke operators.
Jenkins' result recovers the congruences obtained by Edixhoven in the case $m=1$.
More precisely, he proves the following.

\begin{JE}
If $\chi_d(p)=\chi_D(p) \neq 0$, then 
\begin{equation} \label{jen}
A(D,p^{2n}d) = p^nA(p^{2n}D,d),
\end{equation}
and, therefore, 
\textup{(\ref{target})} holds for $m=1$.
\end{JE}

In this paper we prove the following congruences.

\begin{thm} \label{resc}
Let $-d$ and $D$ be fundamental discriminants. 

\textup{ \bf a.} If $\chi_d(p)=\chi_D(p)$, then \textup{(\ref{target})} holds for any $m \geq 1$.

\textup{ \bf b.} If  $\chi_d(p)=-\chi_D(p) \neq 0$, then for any $m \geq 1$
$$
A_m(D,p^{2n+2}d) - A_m(D,p^{2n}d) \equiv 0 \bmod p^n.
$$

\end{thm}

Our argument is also elementary and uses nothing but the identity (\ref{AB}) and some facts about 
the action of half-integral weight Hecke operators. The assumption that both discriminants are
fundamental is not essential. It, however, allows to simplify and generalize the argument.
For instance, Jenkins' identity (\ref{jen}) under this assumption follows at once from the
definitions (see (\ref{sysequ}) below). 
Note that our Theorem \ref{resc}a  essentially 
coincides with the result of Edixhoven. This provides a reason to speculate
that the arithmetic of half-integral weight Hecke operators is somehow connected
with the $p$-adic geometry. Such a connection, if it really exists, looks as an enticing
subject to investigate.

We now turn to the latter part of Ono's question. He noticed \cite[Example 7.15]{Onobook} that, 
for $p \leq 11$, the maximum congruence modulus in (\ref{AO}) exceeds $p$ and called for 
explanations. Recently Boylan \cite{Boylan} found a pretty exact answer in the case $p=2$.
Combining the result of Jenkins, a recent result of the author \cite{Guerzhoy} 
and a theorem of Serre we 
prove the following qualitative result, which uniformly explains the phenomenon without
providing information about the specific congruence moduli.

\begin{thm} \label{smallp}
Let $p \leq 11$. 
If $\chi_d(p)=\chi_D(p) \neq 0$, then the $p$-adic limit 
$$
\lim_{n \rightarrow \infty}  p^{-n} A(D,p^{2n}d) = 0.
$$
\end{thm}



\vspace{3mm}
\begin{center}{\sc 2. Proofs} \end{center}
\vspace{3mm}

Theorem \ref{resc} follows at once from the following proposition.

\begin{prop} \label{id}
Let $p$ be a prime, and let $m \geq 1$ and $n \geq 0$ be integers.
Let $-d$ and $D$ be fundamental discriminants (i.e. 
$-d$ is the discriminant of 
$\Q(\sqrt{-d})$, and $D$ is either $1$ or the discriminant of 
$\Q(\sqrt{D})$).

\textup{ \bf a.} If $\chi_d(p)=\chi_D(p)$ then
\begin{equation} \label{first}
A_m(D,p^{2n}d) = p^nA_m(p^{2n}D,d).
\end{equation}

\textup{ \bf b.} If  $\chi_d(p)=-\chi_D(p) \neq 0$ then
$$
A_m(D,p^{2n+2}d) - A_m(D,p^{2n}d) = p^{n+1}A_m(p^{2n+2}D,d) + p^{n}A_m(p^{2n}D,d).
$$

\end{prop}
 
\begin{proof}[Proof of Proposition \ref{id}]

Let $F=\sum b(n)q^n$ be a (holomorphic or nearly holomorphic) 
modular form of weight $k + 1/2$ for an integer $k \geq 0$, which belongs to the 
Kohnen plus-space (i.e. $b(n)=0$ if $(-1)^k n \equiv 2,3 \bmod 4$). The Hecke operator
$T^+(m)$ acts on 
the Kohnen plus-space and sends $F$ to $F \vert_{k+1/2} T^+(m) = \sum b^+_m(n) q^n$.
If $(-1)^kn$ is a fundamental discriminant, then
\begin{equation} \label{Hecke}
b_m^+(n) = \sum_{l\vert m} \left( \frac{(-1)^kn}{l} \right) l^{k-1} b\left( \frac{m^2}{l^2}n \right).
\end{equation}
This formula follows from \cite[Th. 4.5]{EZ}, where it is proved in the equivalent language of 
Jacobi forms. Although formally the quoted theorem applies only to holomorphic half-integral
weight modular forms, as it is mentioned in \cite{Zagier}, nothing changes if we 
allow the pole at infinity. 
Also note that the technique of Jacobi forms from \cite{EZ} 
covers only the case of odd $k$. However, if $k=0$ (the only even $k$ which we need here) an equivalent
formula may be found in Zagier's paper \cite[proof of Th. 7]{Zagier}. If $k=0$, non-trivial
denominators apparently appear in the right-hand side of (\ref{Hecke}). In order to get rid of these denominators
and to keep our notations compatible with those of \cite{Zagier} we renormalize the Hecke operators by
$$
T(m) = \begin{cases} 
mT^+(m) & \text{if $k=0$} \\
T^+(m) & \text{otherwise}.
\end{cases}
$$
It follows from (\ref{Hecke}) and the definition of the quantities $A_{p^n}(D,d)$ and $B_{p^n}(D,d)$ that
under the assumptions of Proposition \ref{id} 
\begin{align}
A_{p^n}(D,d) = \sum_{i=0}^n\chi_D(p^{n-i})p^iA(p^{2i}D,d)
\notag\\
\label{X} \\
B_{p^n}(D,d) = \sum_{i=0}^n \chi_d(p^{n-i}) B(D,p^{2i}d).
\notag
\end{align}
These equations
combined with (\ref{AB}) imply that for any $n \geq 0$
\begin{equation} \label{sysequ}
\sum_{i=0}^n\chi_D(p^{n-i})p^iA(p^{2i}D,d) = \sum_{i=0}^n\chi_d(p^{n-i})A(D,p^{2i}d),
\end{equation}
and an induction argument in $n$ finishes the proof of Proposition \ref{id} in the case
when $m=1$. 

We now generalize the above argument to the case of arbitrary integer $m \geq 1$.
Recall the usual relation between Hecke operators acting on the 
Kohnen plus-space of
modular forms
of weight $k+1/2$ (see i.e. \cite[Cor. 1 to Th. 4.5]{EZ}):
\begin{equation} \label{relation}
T^+(u)T^+(u')=\sum_{c\vert(u,u')} c^{2k-1} T^+(uu'/c^2).
\end{equation}
In particular, for $k=0,1$ and integers $n, s \geq 0$
$$
\sum_{i=0}^{\min(n,s)} p^iT(p^{n+s-2i}) = T(p^n)T(p^s),
$$
and the equations (\ref{X}) generalize to
\begin{align}
\sum_{i=0}^{\min(n,s)} p^i A_{p^{n+s-2i}}(D,d) = \sum_{i=0}^n\chi_D(p^{n-i})p^iA_{p^s}(p^{2i}D,d)
\notag\\
\label{XX} \\
\sum_{i=0}^{\min(n,s)} p^i B_{p^{n+s-2i}}(D,d) = \sum_{i=0}^n \chi_d(p^{n-i}) B_{p^s}(D,p^{2i}d).
\notag
\end{align}
As previously, (\ref{AB}) implies that the left-hand sides of (\ref{XX}) 
are equal by absolute value and have opposite signs, and
we obtain the following generalization of (\ref{sysequ}) for $n,s \geq 0$
\begin{equation} \label{sysequps}
\sum_{i=0}^n\chi_D(p^{n-i})p^iA_{p^s}(p^{2i}D,d) = \sum_{i=0}^n\chi_d(p^{n-i})A_{p^s}(D,p^{2i}d).
\end{equation}
An induction argument in $n$ finishes the proof of Proposition \ref{id} in the case
when $m=p^s$ with $s>0$. 

Assume now that $m=p^sm_0$ with $p \nmid m_0$ and $s \geq 0$. It follows from 
(\ref{relation}) that
$$
T(p^sm_0)=T(p^s)T(m_0)
$$
for any $k \geq 0$.
It follows that we can multiply the indices in (\ref{sysequps}) and (\ref{XX}) by
$m_0$, as we have had begun the whole argument with the consideration of 
$f_d \vert_\frac{1}{2} T(m_0)$ and $g_D \vert_\frac{3}{2} T(m_0)$ instead of
$f_d$ and $g_D$ correspondingly. This implies the following generalization 
of (\ref{sysequ}),(\ref{sysequps}) for $m_0 \geq 1$, $p \nmid m_0$, $s \geq 0$ and 
any $n \geq 0$:
$$
\sum_{i=0}^n\chi_D(p^{n-i})p^iA_{m_0p^s}(p^{2i}D,d) = \sum_{i=0}^n\chi_d(p^{n-i})A_{m_0p^s}(D,p^{2i}d),
$$
and an induction argument in $n$ completes the proof of Proposition \ref{id}.

\end{proof}


\begin{proof}[Proof of Theorem \ref{smallp}]

It follows from \cite[Theorem 3b]{Guerzhoy} that, since $p \nmid dD$, the formal power series in $q$ 
$$
F = 
\sum_{ \itop{n >0}  {p|n}}\ \  \sum_{\itop{l|n}  {(p,l)=1}} \left( \frac{D}{l} \right) l^{-1} 
A \left(\frac{n^2}{l^2}D,d \right) q^n 
$$
is a $p$-adic cusp form of weight $0$. 
A result of Serre \cite[Th. 7; Rem., p.216]{Serrez} implies 
that, for $p \leq 11$, the $p$-adic limit 
$$
\lim_{n \rightarrow \infty} F \vert U^n = 0,
$$
where $U$ denotes Atkin's $U$-operator $\left( \sum a(n) q^n \right) \vert U = \sum a(pn)q^n$.
Thus the coefficient of $q^{p^n}$ in $F$ approaches $0$ $p$-adically
as $n \rightarrow \infty$.
That is
$$
\lim_{n \rightarrow \infty} A(p^{2n}D,d) =0.
$$
The latter equality combined with Jenkins' identity (\ref{jen}) completes the proof of Theorem \ref{smallp}.
\end{proof}


\begin{thebibliography}{99}
\bibitem{AO}  Ahlgren, Scott; Ono, Ken, Arithmetic of singular moduli and class polynomials, 
Compos. Math. 141 (2005), no. 2, 293-312.
\bibitem{Boylan} Boylan, Matthew, $2$-adic properties of Hecke traces of singular moduli,
Math. Res. Lett., 12 (2005), 10001-10014. 
\bibitem{Edixhoven} Edixhoven, Bas, On the $p$-adic geometry of traces of singular moduli,
preprint, 2005, arxiv math.NT/0502213 v1.
\bibitem{EZ}  Eichler, Martin; Zagier, Don, The theory of Jacobi forms, 
Progress in Mathematics, 55. Birkh\"auser Boston, Inc., Boston, MA, 1985.
\bibitem{Guerzhoy} Guerzhoy, Pavel, The Borcherds-Zagier isomorphism and a $p$-adic version of
the Kohnen-Shimura map,  Int. Math. Res. Not. 2005, no. 13, 799-814.
\bibitem{Jenkins} Jenkins, Paul, $p$-adic properties for traces
of singular moduli, accepted for publication in the
  International Journal of Number Theory
\bibitem{Onobook} Ono, Ken, The Web of Modularity: Arithmetics of the Coefficients
of Modular Forms and $q$-series, 
CBMS Reg. Conf. Ser. Math., 102, 
Amer. Math. Soc. Providence RI, 2004.
\bibitem{Serrez} Serre, Jean-Pierre, 
Formes modulaires et fonctions z\^eta $p$-adiques,  
Modular Functions of One Variable, III, Proc. Internat. Summer School, Univ. Antwerp, 1972, pp. 191--268, Lecture
Notes in Math., Vol. 350, 
Springer-Verlag, Berlin, 1973.
\bibitem{Zagier} Zagier, Don, Traces of singular moduli, 
Motives, Polylogarithms and Hodge Theory, Part I (Irvine, CA, 1998), 
211--244, Int. Press Lect. Ser., 3, I, Int. Press, Somerville, MA, 2002
\end{thebibliography}
\end{document}